\documentclass[12pt, reqno]{amsart}
\usepackage{amsmath, amsthm, amscd, amsfonts, amssymb, graphicx, color}
\usepackage[bookmarksnumbered, colorlinks, plainpages]{hyperref}

\newtheorem{theorem}{Theorem}[section]
\newtheorem{lemma}[theorem]{Lemma}
\newtheorem{proposition}[theorem]{Proposition}
\newtheorem{corollary}[theorem]{Corollary}
\theoremstyle{definition}
\newtheorem{definition}[theorem]{Definition}
\newtheorem{example}[theorem]{Example}

\theoremstyle{remark}
\newtheorem{remark}[theorem]{Remark}
\numberwithin{equation}{section}

\begin{document}
\setcounter{page}{1}

\title[$q$-Frequent hypercyclicity]{$q$-Frequently hypercyclic operators}

\author[M. Gupta and A. Mundayadan]{Manjul Gupta$^{*}$ and Aneesh Mundayadan}

\address{$^{1}$ Department of Mathematics and Statistics, Indian Institute of Technology Kanpur, Kanpur 208 016, UP, India.}
\email{\textcolor[rgb]{0.00,0.00,0.84}{manjul@iitk.ac.in}}
\email{\textcolor[rgb]{0.00,0.00,0.84}{aneeshm@iitk.ac.in;aneeshkolappa@gmail.com}}




\date{Received: xxxxxx; Revised: yyyyyy; Accepted: zzzzzz.
\newline \indent $^{*}$ Corresponding author}

\begin{abstract}
We introduce $q$-frequently hypercyclic operators and derive a
sufficient criterion for a continuous operator to be $q$-frequently
hypercyclic on a locally convex space. Applications are given to
obtain $q$-frequently hypercyclic operators with respect to the
norm-, $F$-norm- and weak*- topologies. Finally, the frequent
hypercyclicity of the non-convolution operator $T_\mu$ defined by
$T_\mu(f)(z)=f'(\mu z)$, $|\mu|\geq 1$ on the space $H(\mathbb{C})$
of entire functions equipped with the compact-open topology is
shown.
\end{abstract} \maketitle

\section{Introduction}

The main theme in the dynamics of linear operators is the notion of
hypercyclicity which plays an important role in the study of the
invariant subset problem in Banach spaces. This notion was initiated
by S. Rolewicz \cite{Rolewicz} in the setting  of infinite
dimensional Banach spaces in 1969, though the examples of
translation and differential operators on the space of entire
functions equipped with the compact-open topology were known to be
hypercyclic in an earlier work of G. D. Birkhoff \cite{Birkhoff} and
G. R. MacLane \cite{MacLane}. Now a vast literature dealing with
hypercyclicity of  operators as well as other related notions in
linear dynamics is available in \cite{Bayart Matheron},
\cite{Erdmann Peris2}, and \cite{Erdmann}.

In 2006, F. Bayart and S. Grivaux \cite{Bayart Grivaux} further
strengthened this concept to frequent hypercyclicity, which
quantifies the frequency with which the iterates of a given linear
operator at a point visit each non-empty open set. After the
appearance of this work, several results on frequently hypercyclic
operators have been established, for instance one may refer to
\cite{Blasco Bonilla Erdmann}, \cite{Bonilla Erdmann1},
\cite{Bonilla Erdmann2}, \cite{Erdmann Peris1} and \cite{Shkarin}.
In this paper we introduce $q$-frequent hypercyclicity which lies
between hypercyclicity and frequent hypercyclicity, where $q$ is a
fixed natural number. The case $q=1$ coincides with frequent
hypercyclicity. We prove a sufficient criterion for a continuous
linear operator to be $q$-frequently hypercyclic on a locally convex
space and give applications to obtain $q$-frequently hypercyclic
operators with respect to the norm on Banach spaces, the $F$-norm on
$F$-spaces and the weak*-topology on dual of Banach spaces. We also
provide examples of hypercyclic operators that are not
$q$-frequently hypercyclic for any $q\in \mathbb{N}$.

\section{Preliminaries}

Let $X$ be a separable topological vector space, and
$\mathcal{L}(X)$ denote the space of all continuous linear operators
on $X$. An operator $T\in \mathcal{L}(X)$ is said to be
\emph{hypercyclic} if there exists a vector $x\in X$ such that the
orbit $\{T^n(x):n\geq 0\}$ is dense in $X$. Such a vector $x$ is
called a \emph{hypercyclic vector} for $T$. As mentioned in the
previous section, Birkhoff's translation operator
$T_a(f)(z)=f(z+a)$, for nonzero $a\in \mathbb{C}$  and MacLane's
differentiation operator $D(f)=f'$ on the space $H(\mathbb{C})$ are
hypercyclic. Also, Rolewicz proved the hypercyclicity of  the
operator $\lambda B$ on $\ell^p$ or $c_0$ for $1\leq p <\infty$ and
$|\lambda|
> 1$, where $B$ is the unweighted backward shift defined by
$B(e_n)=e_{n-1}, n\geq1$, with $e_0=0$ and $e_n=\{0,0,..,1,0..\}$,
$1$ being placed at the $n$th coordinate. Generalizing this result,
H.N. Salas \cite{Salas} proved that the weighted shift $B_w$
associated to a weight sequence $(w_n)$ of positive reals, given by
$B_w(e_n)=w_ne_{n-1}, n\geq1$ is hypercyclic on $\ell^p$ or $c_0$ if
and only if $\limsup_{n\to \infty}(w_1w_2...w_n)=\infty$.
\\

For testing the hypercyclicity of a linear operator, a sufficient
criterion known as the hypercyclicity criterion, initially obtained
by Kitai \cite{Kitai}, has appeared in different forms and the one
which is given below is due to H. Petterson \cite{Petterson}. This
is useful even for linear operators defined on non-metrizable
topological vector spaces. For the definition of $F$-norm, we refer
to \cite{Bonilla Erdmann2}, p. 385.

\begin{theorem} [Hypercyclicity criterion] \label{thm2.1} Let $(X,\tau)$ be
a separable topological vector space. Suppose further that $X$
carries an $F$-norm $||.||$ with respect to which it is complete and
that $\|.\|$-topology is stronger than $\tau$. If $T$ is an operator
continuous with respect to the $F$-norm, $D\subset X$ a countable
$\tau$-dense set and $S_n:D\rightarrow X$ maps such that, for all
$x\in D$,
\begin{itemize}
\item[1.]$||T^{n}(x)||\rightarrow 0$ and $||S_n(x)||\rightarrow 0$ as $n\rightarrow \infty$; and
\item[2.]$T^nS_n(x)=x$, for each $n\in \mathbb{N}$,
\end{itemize}
then the operator $T$ is  $\tau$-hypercyclic.
\end{theorem}

An operator $T\in \mathcal{L}(X)$ is called \emph{frequently
hypercyclic} if there exists an $x$, called a \emph{frequently
hypercyclic vector} for $T$ such that for every nonempty open set
$U$ in $X$, the set $\mathbf{N}(x,U)=\{ n \in \mathbb{N}: T^n(x) \in
U\}$ has positive lower density; where the \emph{lower density} of
a subset $A$ of $\mathbb{N}$, the set of natural numbers, is defined
as
\begin{equation}
\begin{gathered}
\underline{\text{dens}}(A) $ = $\displaystyle \liminf_{N\to \infty}
\dfrac{\text{card}\{n\in A: n\leq N\}}{N},
\end{gathered}
\end{equation}
the symbol $\text{card}(B)$ being used to denote the cardinality of
the set $B$.

Let us note that $\underline{\text{dens}}(B)\in [0,1]$ for any
subset $B$ of $\mathbb{N}$. Clearly, the lower density of any finite
set is zero and that of $\mathbb{N}$ is 1. If $A$ is a strictly
increasing sequence $(n_k)$, the lower density of $A$ is
characterized as \cite{Bonilla Erdmann2}
\begin{equation}
\underline{\text{dens}}(n_k)  =  \liminf_{k\to \infty}
\frac{k}{n_k}.
\end{equation}
Alternatively, an operator $T\in L(X)$ is frequently hypercyclic if
there is some $x\in X$ such that for every nonempty open subset $U$
of $X$, there exist a strictly increasing $(n_k)$ of natural numbers
and a constant $C>0$ such that
\begin{center}
$T^{n_k}(x)\in U$ $\text{and}$ $n_k \leq Ck,$ $\forall$ $k\in
\mathbb{N}.$
\end{center}

Analogous to the hypercyclicity criterion, we have the following
criterion, proved in \cite{Bayart Grivaux} and \cite{Bonilla
Erdmann2}.
\begin{theorem}[Frequent hypercyclicity criterion]\label{thm2.2}
Let $X$ be a separable $F$-space and $T\in L(X)$. If there exist a
dense subset $D\subset X$ and a map $S:D\rightarrow D$ such that
\begin{itemize}
\item[1.]$\sum T^n(x)$ and $\sum S^n(x)$ are unconditionally convergent
for each $x\in D$; and
\item[2.] $TS=I$, the identity on $D$,
\end{itemize}
then the operator $T$ is frequently hypercyclic.
\end{theorem}

The operators of Birkhoff and MacLane satisfy the above criterion,
cf. \cite{Bayart Grivaux} and \cite{Bayart Matheron}, and so they
are frequently hypercyclic. In fact, any continuous operator, except
a scalar multiple of the identity, that commutes with all
translations on $H(\mathbb{C})$, has been shown to be frequently
hypercyclic \cite {Bonilla Erdmann1}. We also recall the hypercyclic
comparison principle from \cite{Bayart Grivaux} and \cite{Bayart
Matheron}, which says how to transfer the hypercyclicity via a
linear quasi-conjugacy.

\begin{proposition}[Hypercyclic comparison principle]\label{prop2.3} Let $T$ and $S$
be continuous linear operators on two topological vector spaces $X$
and $Y$ respectively and $A:X\rightarrow Y$ be a continuous linear
map with dense range such that $SA=AT$. If $T$ is hypercyclic (or
frequently hypercyclic) on $X$, then $S$ is hypercyclic (or
frequently hypercyclic) on $Y$.
\end{proposition}

\section {$q$-Frequently hypercyclic operators
} We first introduce the $q$-lower density of a subset of natural
numbers, for $q\in \mathbb{N}$ and determine a useful
characterization.

\begin{definition}
Let $A \subset \mathbb{N}$ and $q\in \mathbb{N}$. The $q$-lower
density of $A$ is defined as
\begin{center}
$q$-$\underline{\text{dens}}(A) = \displaystyle \liminf_{N\to
\infty} \dfrac{\text{card}\{n\in A: n\leq N^q \}}{N}.$
\end{center}
\end{definition}

Let us note that the lower density of a set is always finite and
lies in the interval $[0,1]$, but the $q$-lower density can vary in
$[0,\infty]$ for $q\geq 2$. As in the case of the lower density, we
have the following.

\begin{proposition}\label{lm2}
Let $(n_k)$ be a strictly increasing sequence of natural numbers.
Then
\begin{enumerate}
\item $q$-$ \emph{\underline{dens}}(n_k) = \displaystyle \liminf_{k\to
\infty} \dfrac {k}{{n_k}^{1/q}}$.
\item $q$-$\emph{\underline{dens}}(n_k)>0$ if and only if there exists a
constant $C>0$ such that $n_k \leq Ck^q$ for all $k\in \mathbb{N}$.
\end{enumerate}
\end{proposition}

\begin{proof}
Fix a number $k\in \mathbb{N}$. For any $N\in \mathbb{N}$ with
$n_k\leq N^q<n_{k+1}$, we have that
\begin{center}
$p_N=\dfrac{\text{card}\{k\in \mathbb{N}: n_k\leq N^q \}}{N}=k/N$.
\end{center}
 Thus the inequality
\begin{center}
$\dfrac {k}{{n_{k+1}}^{1/q}}<  p_N \leq \dfrac {k}{{n_k}^{1/q}}$
\end{center}
implies the first part in the theorem.

Part $(2)$ follows immediately from the fact that $\liminf a_k>0$ if
and only if $\frac {1} {a_k} \leq C$ for some $C>0$, where $(a_k)$
is a sequence of positive numbers.
\end{proof}

We now define the notion of $q$-frequent hypercyclicity of linear
operators on topological vector spaces.
\begin{definition}
Let $q \in  \mathbb{N}$. A continuous linear operator $T$ on a
separable topological vector space $X$ is said to be $q$-frequently
hypercyclic if there exists an $x \in X$ such that for any nonempty
open subset $U$ of $X$, the set $\mathbf{N}(x,U)=\{ n \in
\mathbb{N}: T^n(x) \in U\}$ has positive $q$-lower density. Such a
vector is called a $q$-frequently hypercyclic vector for $T$.
\end{definition}

Alternatively, a continuous linear operator $T$ on a separable
topological vector space $X$ is \emph{$q$-frequently hypercyclic}
if there exists an $x \in X$ such that for any nonempty open subset
$U$ of $X$, we can find a strictly increasing sequence $(n_k)$ of
natural numbers and a constant $C>0$ such that

\begin{center}
$T^{n_k}(x)\in U$ and $n_k \leq Ck^q,$ for all $k\in \mathbb{N}.$
\end{center}
Such a vector is called a \emph{$q$-frequently hypercyclic vector}
for $T$.

Obviously, every frequently hypercyclic operator is $q$-frequently
hypercyclic for any natural number $q$, and the two notions are the
same for the case $q=1$. Also this new property of linear operators
is stronger than hypercyclicity; however, none of the converse
implications is true, e.g. consider

\begin{example}
Here we show that there exists a hypercyclic operator on $\ell^1$
that is not 2-frequently hypercyclic with respect to the weak
topology of $\ell^1$. Indeed, the weighted backward shift $B_w$ with
weights $w_n=\sqrt{\frac{n+1}{n}}$ is hypercyclic by the result of
Salas, but was shown to be non-frequently hypercyclic on $\ell^2$ in
\cite{Bayart Grivaux}. For showing the non-2-frequent hypercyclicity
of $B_w$, choose the weakly open set $U=\{(y_n)\in
\ell^1:|y_1|>1\}$. Let $x=(x_n)$ be a $2$-frequently hypercyclic
vector for $B_w$. Enumerate the set $N(x,U)=\{ n \in \mathbb{N}:
B_w^n(x) \in U\}$ as $(n_k)$. Thus we have a constant $c>0$ such
that $n_k \leq ck^2$, and hence

\begin{center}
$\displaystyle \sum_{k\geq 1}\frac {1}{\sqrt{n_k}}=\infty$.
\end{center}
On the other hand, $B_w^{n_k}(x) \in U$ implies that

\begin{center}
 $ \sqrt{n_k+1}$ $ |x_{n_k+1}|> 1$, for all $k\geq 1$
\end{center}As $(x_n)\in \ell^1$, we get

\begin{center}
$\sum_{k\geq 1}\frac {1}{\sqrt{n_k+1}} <\infty$,
\end{center}
which is a contradiction. Hence $B_w$ is not 2-frequently
hypercyclic on $\ell^1$ for the norm topology.
\end{example}

\begin{example}\label{example3.4} We show the existence of a hypercyclic operator that is not
$q$-frequently hypercyclic with respect to the weak topology, for
any $q\in \mathbb{N}$. Let us consider the unilateral shift $B_w$ on
$\ell^1$ with weights $w_k=\displaystyle \frac
{\ln(k+2)}{\ln(k+1)}$, $k\in \mathbb{N}$. By the result of
H.N.Salas, $B_w$ is hypercyclic since $w_1w_2...w_k= \frac
{\ln(k+2)}{\ln 2} \rightarrow \infty$ as $k \rightarrow \infty$. Let
$x=(x_n)$ be a $q$-frequently hypercyclic vector for $B_w$ for some
$q\in \mathbb{N}$. Enumerate the set $N(x,U)=\{ n \in \mathbb{N}:
B_w^n(x) \in U\}$ as $(n_k)$, where $U=\{(y_n)\in \ell^1:|y_1|>1\}$.
Thus we have a constant $c>0$ such that $n_k \leq ck^q$. This
implies that $\ln(n_k)\leq C \ln k$,
 for some constant $C>0$ and for all $k\in \mathbb{N}$. Hence

\begin{center}
$\displaystyle \sum_{k\geq 1}\frac {1}{\ln{n_k}}=\infty$.
\end{center}
On the other hand, $B_w^{n_k}(x) \in U$ implies that

\begin{center}
$\displaystyle \frac {\ln({n_k}+2)}{\ln 2} |x_{n_k+1}|> 1$.
\end{center} Consequently,
\begin{center}
$\displaystyle \sum_{k\geq 1}\frac {1}{\ln{(n_k+2)}} <\infty$,
\end{center} which is a contradiction. Hence $B_w$ is not
$q$-frequently hypercyclic in the norm topology, for any $q\in
\mathbb{N}$.
\end{example}

We now prove a criterion, similar to the frequent hypercyclicity
criterion, which works even for operators defined on certain
non-metrizable locally convex spaces. Using this, we obtain a
2-frequently hypercyclic operator that is not frequently
hypercyclic. Before stating the result, let us recall that a series
$\sum {x_n}$ in a topological vector space is
\emph{unconditionally convergent} if $\sum {x_{\sigma (n)}}$ is
convergent for every permutation $\sigma$ of $\mathbb{N}$. In any
topological vector space, this mode of convergence is equivalent to
the unordered convergence of $\sum {x_n}$, cf. \cite{Kamthan Gupta},
p.154. Thus a series $\sum {x_n}$ is unconditionally convergent if
and only if for every non-empty open set $U$ of $0$, there
corresponds an $N\in \mathbb{N}$ such that $\sum _{n\in F}{x_n}\in
U$ for every finite set $F\subset [N,\infty)$. Also, for the proof
of our criterion we need the following lemma from \cite{Bayart
Matheron}.

\begin{lemma} \label{lem3.4} Let $(N_k)$ be a strictly increasing sequence of
natural numbers. Then there exists a pairwise disjoint sequence
$(J_k)$ of subsets of $\mathbb{N}$ such that
\begin{enumerate}
\item $\emph{\underline{dens}}(J_k)>0$ for each $k\geq 1$
\item $|n-m|\geq N_k+N_p$ for $n\not=m$ and $(n,m)\in
J_k\times J_p$.
\item $n\geq N_k$, for each $n\in J_k$ and $k\geq 1$.
\end{enumerate}
\end{lemma}

We now state and prove the main result of the paper.

\begin{theorem}[$q$-frequent hypercyclicity criterion]
\label{thm3.5} Let $(X,\tau)$ be a separable locally convex space
and $q\in \mathbb{N}$. Suppose that $X$ is equipped with  an
$F$-norm $\|.\|$ such that the F-norm topology is stronger than
$\tau$ and $(X,\|.\|)$ is complete. If $T$ is an operator continuous
with respect to the $F$-norm, $D$ is a subset of $X$ containing a
countable $\tau$-dense subset and $S:D\rightarrow D$ is a map such
that

\begin{enumerate}
\item $\sum T^{n^q}(x)$ and $\sum S^{n^q}(x)$ are unconditionally
convergent with respect to the $F$-norm, for each $x\in D$; and
\item $TS=I$, the identity on $D$,
\end{enumerate}
then the operator $T$ is $q$-frequently hypercyclic with respect to
$\tau$.
\end{theorem}

\begin{proof}

Our proof is inspired by that of the frequent hypercyclicity
criterion given in \cite{Bayart Matheron}. However, we outline the
proof for the sake of completeness. Let $Y=\{x_1,x_2,... \}\subset
D$ be a countable $\tau$-dense set. Consider a summable sequence
$(\epsilon _k)$ of positive real numbers, which are to be chosen
later. By the hypothesis, corresponding to $\epsilon_k$, we can find
$N_k\in \mathbb{N}$ such that

\begin{equation}\label{eq1}
\left\|\sum_{n\in F} T^{n^q}(x_i)\right\|+ \left\| \sum_{n\in F}
S^{n^q}(x_i)\right\|<\epsilon_k,  1\leq i\leq k,
\end{equation}
for any finite set $F\subset [N_k,\infty)$ of natural numbers. We
may now assume that $(N_k)$ is strictly increasing so that by Lemma
\ref{lem3.4}, we get a sequence $(J_k)$ of subsets of $\mathbb{N}$
with the properties mentioned therein. We now set

\begin{center}
$\displaystyle x=\sum_{k\geq 1}\sum_{n\in J_k} S^{n^q}(x_k)$.
\end{center}
Since unconditional convergence implies subseries convergence
\cite{Kamthan Gupta} p.154, the series $\sum_{n\in J_k}
S^{n^q}(x_k)$ converges for each natural number $k$. It follows by
\eqref{eq1} that
\begin{center}
 $\displaystyle \sum_{k=1}^{\infty}\left\|\sum_{n\in J_k} S^{n^q}(x_k)\right\| \leq
 \sum_{k\geq 1}
{\epsilon_k}< \infty$
\end{center}
Thus $x\in X$.

Let us now fix $k\in \mathbb{N}$ and $m\in J_k$. Then

\begin{equation*}
\displaystyle T^{m^q}(x)=\sum_{l\geq 1}\sum_{n\in J_l}
T^{m^q}S^{n^q}(x_l).
\end{equation*} and so
\begin{equation}\label{eq2}
\displaystyle \left\|T^{m^q}(x)-x_k\right\|\leq \sum_{l\geq
1}\left\|\sum_{n\in J_l,m>n} T^{m^q-n^q}(x_l)\right\|+\displaystyle
\sum_{l\geq 1}\left\|\sum_{n\in J_l,m<n} S^{n^q-m^q}(x_l)\right\|.
\end{equation}

Let us now consider the first sum on the right hand side of
\eqref{eq2}. Indeed, writing the first term as

\begin{equation*}
\begin{gathered}
\displaystyle \sum_{l=1}^{k}\left\|\sum_{n\in J_l,m>n}
T^{m^q-n^q}(x_l)\right\| +\displaystyle \sum_{l\geq
k+1}\left\|\sum_{n\in J_l,m>n} T^{m^q-n^q}(x_l)\right\|,
\end{gathered}
\end{equation*}
 we have
\begin{equation*}
\begin{gathered}
\displaystyle \sum_{l\geq 1}\left\|\sum_{n\in J_l,m>n}
T^{m^q-n^q}(x_l)\right\| \leq k\epsilon_k +\sum_{j\geq k+1}
\epsilon_j.
\end{gathered}
\end{equation*}
The last inequality arrives because whenever $m\in J_k$, we have
$m^q-n^q>$max$(N_k,N_l)$ for any $n\in J_l$ with $m>n$. Similarly,
we evaluate the second term to get

\begin{equation*}
\displaystyle \sum_{l\geq 1}\left\|\sum_{n\in J_l,m<n}
S^{n^q-m^q}(x_l)\right\| \leq k\epsilon_k +\sum_{j\geq k+1}
\epsilon_j.
\end{equation*} Set $\alpha_k=k\epsilon_k +\sum_{j\geq
k+1}\epsilon_j.$ So we arrive at the inequality,

\begin{equation}\label {eq3}
\displaystyle \left\|T^{m^q}(x)-x_k\right\|\leq 2\alpha_k<3
\alpha_k, \forall m\in J_k, k\geq 1.
\end{equation}
Choose $\epsilon_k$ such that $\alpha_k \rightarrow 0$. We now show
that $x$ is a $q$-frequently hypercyclic vector for the operator $T$
with respect to the topology $\tau$. Let $G$ be a nonempty
$\tau$-open set and $y+U\subset G$ for some $\tau$-neighborhood $U$
of the origin. Then we find a balanced neighborhood $V$ of the
origin such that $V+V\subset U$. Since $Y$ is $\tau$-dense in $X$,
we find an increasing sequence of natural numbers $(n_k)$ such that
$x_{n_k}-y\in V$ for all $k\geq 1$. Since the $||.||$-topology is
finer than $\tau$, it follows from \eqref{eq3} that, for some $N\in
\mathbb{N}$, $T^{m^q}(x)-x_k \in V$ for every $m\in J_k$ and $k\geq
N$. Thus from the facts that $x_{n_N}-y\in V$ and
$T^{m^q}(x)-x_{n_N} \in V$, we obtain that for all $m\in J_{n_N}$,
$T^{m^q}(x)-y \in V+V \subset U$. Our conclusion now follows because
$T^{m^q}(x)\in G $ for all $m\in J_{n_N}$, which has positive lower
density.
\end{proof}

\begin{remark}
It is evident from the proof of Theorem \ref{thm3.5} that the
sequence $(T^{n^q})$ is frequently hypercyclic and thus $T$ is
$q$-frequently hypercyclic. It would be interesting to know whether
the converse is true or not, i.e., is $(T^{n^q})$ frequently
hypercyclic whenever $T$ is $q$-frequently hypercyclic?
\end{remark}

\begin{remark} \label{rem3.6} Let us also
note that in the above theorem, the countability assumption on a
subset of $D$ may be waived in case the topology $\tau$ is generated
by an $F$-norm. Indeed, if $(x_n)$ is a $\tau$-dense sequence in
$X$, choose a countable set $\{y_{n,m}:n,m\geq1\}$ where $y_{n,m}\in
D\cap B(x_n,\frac {1} {m})$,  $B(x_n,\frac {1} {m})$ being the open
ball of radius $\frac {1} {m}$ centered at $x_n$. It is now easy to
verify that $\{y_{n,m}:n,m\geq1\}$ is a $\tau$-dense set in $X$.
\end{remark}

\section{Applications}

In this section, we consider some applications of the $q$-frequent
hypercyclicity criterion for obtaining $q$-frequently hypercyclic
operators on spaces equipped with linear topologies which are not
necessarily metrizable. Besides, we prove the frequent
hypercyclicity of a non-convolution operator on $H(\mathbb{C})$, at
the end of the section. Let us begin with the results on sequence
spaces with metrizable topologies.

\begin{proposition}\label{prop4.1}
Let $\lambda$ be a sequence space equipped with an $F$-norm and let
$\{e_n\}$ be an unconditional basis in $\lambda$. If for some $q\in
\mathbb{N}$,
\begin{center}
$\displaystyle \sum_{n\in\mathbb{N}} \frac {1}
{w_1w_2...w_{n^q+j}}e_{n^q+j}$
\end{center}
converges unconditionally for each $ j\in \mathbb{N}$, then the
backward shift $B_w$ associated to the weight sequence $(w_n)$, is
$q$-frequently hypercyclic.
\end{proposition}

\begin{proof}
Since $\lambda$ is an $F$-space, in view of Remark \ref {rem3.6}
choose $D$ to be the dense set spanned by $\{e_n:n\geq1\}$. Define
$S_w$ on $D$ by $S_w(e_n)= \frac {1} {w_{n+1}} e_{n+1}$. To apply
our criterion, we are only required to prove the unconditional
convergence of $\sum S_w^{n^q}(x)$ for each $x\in D$. Indeed, for a
given $k\in \mathbb{N}$, we have
\begin{center}
$\displaystyle S_w^n(e_k)= \frac {1} {w_{k+1}...w_{k+n}} e_{k+n}$.
\end{center}
Thus, by the hypothesis, the series
\begin{center}
$\displaystyle \sum _{n\geq 1} S_w^{n^q}(e_k)=w_1w_2...w_k
\sum_{n\in\mathbb{N}} \frac {1} {w_1...w_{k+n^q}} e_{k+n^q}$.
\end{center}
converges unconditionally in $\lambda$. Hence $B_w$ is
$q$-frequently hypercyclic.

\end{proof}

As a consequence of the above result, we obtain a $q$-frequently
hypercyclic operator that is not frequently hypercyclic. This is the
Bergman shift, considered in \cite {Bayart Grivaux}.

\begin{corollary}
Let $B_w$ be the unilateral shift on $\ell^2$ given by the weights
$w_n=\sqrt{\frac{n+1}{n}}$, $n\geq 1$. Then $B_w$ is $2$-frequently
hypercyclic and is not frequently hypercyclic; a fortiori, $B_w$ is
$q$-frequently hypercyclic for any $q\geq 2$.
\end{corollary}
\begin{proof}
Since $w_1w_2...w_{n^2+j}=\sqrt{n^2+j+1}$, the result follows.
\end{proof}

\begin{remark}
One may apply Proposition \ref{prop4.1} to obtain the $q$-frequent
hypercyclicity of shift operators defined on Fr\'{e}chet spaces, for
example the space $H(\mathbb{C})$ of entire functions equipped with
the compact-open topology, the space of all sequences with the
topology of co-ordinate convergence and the classical $\ell^p$
spaces.
\end{remark}

Before we move on to another application, let us see an example.
\begin{example} Let $p\in \mathbb{N}$. Then there
exists an operator which is not $p$-frequently hypercyclic, but it
is $q$-frequently hypercyclic for all $q\geq p+1$. We consider the
unilateral shift $B_w$ on $\ell^2$ with weights $w_k=\displaystyle
{\left(\frac {k+2}{k+1}\right)}^{1/2p}$, $k\in \mathbb{N}$. Then the
proof similar to that of Example \ref{example3.4} shows that $B_w$
is not $p$-frequently hypercyclic on $\ell^2$. We can also conclude
that $B_w$ is $(p+1)$-frequently hypercyclic by applying Proposition
\ref{prop4.1} to the dense set of finite sequences.
\end{example}

Our next application of Theorem \ref{thm3.5} is for the bilateral
backward shift operators. Let $X$ be an $F$-sequence space over the
set $\mathbb{Z}$ of integers such that the unit sequences
$(e_n)_{n\in \mathbb{Z}}$ form an unconditional basis in $X$. For
$w=(w_n)\in \ell^\infty(\mathbb{Z})$, the operator
$T_w(e_n)=w_ne_{n-1}$ is the bilateral backward shift. The content
of the following proposition is the $q$-frequent hypercyclicity of
$T_w$.

\begin{proposition}\label{prop4.5}
Let $(e_n)_{n\in \mathbb{Z}}$ form a unconditional basis in an
$F$-sequence space $X$ and let $q\in \mathbb{N}$. If
\begin{center}
$\displaystyle \sum_{n\in\mathbb{N}} \frac {1}
{w_1w_2...w_{n^q+j}}e_{n^q+j}$ and $\displaystyle
\sum_{n\in\mathbb{N}} w_jw_{j-1}...w_{-n^q+j+1}e_{-n^q+j}$
\end{center}
converge unconditionally for each $ j\in \mathbb{N}$, then $T_w$ is
$q$-frequently hypercyclic on $X$.

\end{proposition}

\begin{proof}
Let $D$ be the set spanned by the sequence $(e_n)_{n\in\mathbb{Z}}$.
Consider the map $S_w(e_k)=\frac {1} {w_{k+1}}e_{k+1}$ on the dense
set $D$ of $X$, so that $T_wS_w$ is the identity operator on $D$.
Then,
\begin{center}
$T_w^{n^q}(e_j)=w_jw_{j-1}..w_{j-n^q+1}e_{j-n^q}$
\end{center} and
\begin{center}
$S_w^{n^q}(e_j)=\frac {1} {w_{j+1}w_{j+2}...w_{j+n^q}}e_{j+n^q}$
\end{center} for each $j\in \mathbb{Z}$.
From the hypothesis, we obtain that the series $\sum T_w^{n^q}(e_j)$
and $\sum S_w^{n^q}(e_j)$ converge unconditionally in $X$. The
desired result now follows since $T_w$ and $S_w$ are linear and $D$
is the span of $(e_n)_{n\in\mathbb{Z}}$.
\end{proof}

A particular case of Proposition \ref{prop4.5} is when $T_w$ is the
bilateral backward shift defined on the sequence space
$\ell^p(\mathbb{Z})$ or $c_0(\mathbb{Z})$ for $1\leq p<\infty$.
Recall that hypercyclicity of $T_w$ was characterized by H. N. Salas
in \cite{Salas}. Also a series $\sum_{n\in \mathbb{Z}} a_ne_n$
converges unconditionally in $\ell^p(\mathbb{Z})$ if and only if the
sequence $(a_n)\in \ell^p(\mathbb{Z})$. Thus we derive the following
result from Proposition \ref{prop4.5}.

\begin{corollary}\label{cor4.6}
Let $q\in \mathbb{N}$ and $1\leq p<\infty$. Assume that for each
$j\in \mathbb{Z}$,
\begin{center}
$\sum_{n\in\mathbb{N}} \frac {1} {(w_1w_2...w_{n^q+j})^p}<\infty$
and $\displaystyle \sum_{n\in\mathbb{N}}
(w_jw_{j-1}...w_{-n^q+j+1})^p<\infty$.
\end{center}
Then $T_w$ is $q$-frequently hypercyclic on $\ell^p(\mathbb{Z})$. If
$\lim_{n\rightarrow \infty} (w_1w_2...w_{n^q+j})= \infty$ and
$\lim_{n\rightarrow \infty} (w_jw_{j-1}...w_{-n^q+j+1})=0$ for each
$j\in \mathbb{Z}$, then $T_w$ is $q$-frequently hypercyclic on
$c_0(\mathbb{Z})$.
\end{corollary}

We have yet another application of Theorem \ref{thm3.5}. Let $X$ be
a Banach space having a Schauder basis $\{x_n,f_n\}$. Then the dual
$X^*$ is weak*-separable and the weak*-topology is not metrizable,
when $X$ is infinite dimensional. We assume that $\{x_n,f_n\}$ is a
\emph{symmetric Schauder basis} for $X$, which means that $ \sum
_{n\geq 1} f_{\mu(n)}(x)x_{\sigma(n)}$ converges for each $x\in X$
and every pair $(\mu,\sigma)$ of permutations of $\mathbb{N}$. A
symmetric base in a Banach space is regular $(\inf ||x_n||>0)$ and
bounded $(\sup ||x_n||< \infty)$, see for example \cite {Kamthan
Gupta1}, p. 133. Corresponding to a weight sequence $(w_n)$ and a
symmetric Schauder basis $\{x_n,f_n\}$ (where the index starts from
1), we define the backward shift $B_w$ by $B_w(f_n)=w_nf_{n-1}$,
$n\geq1$ with $f_0=0$. This operator is continuous with respect to
the norm as well as the weak*-topology of $X^*$, cf. \cite
{Petterson}. We now prove:

\begin{theorem} \label{thm4.4}
Let $X$ be a Banach space with a symmetric Schauder basis
$\{x_n,f_n\}$. Then for $q\in \mathbb{N}$, the backward shift
operator on $X^*$ is $q$-frequently hypercyclic with respect to the
weak*-topology of $X^*$ if $\displaystyle \sum_{n\in\mathbb{N}}
\frac {1} {w_{1}...w_{j+n^q}}$ converges for each $j\in \mathbb{N}$.
In particular, if $\displaystyle \sum_{n\in\mathbb{N}} \frac {1}
{w_{1}w_{2}...w_{n}}$ converges, then $B_w$ is weak*-frequently
hypercyclic on $X^*$.
\end{theorem}

\begin{proof}
In order to apply Theorem \ref{thm3.5}, consider $D$ to be the span
of $\{f_n:n\geq 1\}$. Then $D$ contains a countable weak$^*$-dense
subset. The forward shift $S_w(f_n)= \frac {1} {w_{n+1}} f_{n+1}$,
$n\geq 1$ maps $D$ to itself. Thus
\begin{center}
$\displaystyle S_w^n(f_k)= \frac {1} {w_{k+1}...w_{k+n}} f_{k+n}$
\end{center} and
\begin{center}
$\displaystyle \sum_{n\in \mathbb{N}} S_w^{n^q}(f_k)=w_1w_2...w_k
\sum_{n\in\mathbb{N}} \frac {1} {w_{1}...w_{k+n^q}} f_{k+n^q}$.
\end{center}

Since a symmetric Schauder basis is regular, we have that
$||f_n||<K$ for some constant $K>0$ and for all $n\in \mathbb{N}$,
cf. \cite{Kamthan Gupta1}, p. 261 and \cite{Singer}, p. 25. Hence by
our hypothesis, the series $\sum_{n\in \mathbb{N}} S_w^{n^q}(f_k)$
is absolutely convergent and so unconditionally convergent.
Consequently, the shift $B_w$ is $q$-frequently hypercyclic with
respect to the weak*-topology on $X^*$ by Theorem \ref{thm3.5}.
\end{proof}

As a consequence of the above result, we derive:

\begin{corollary}\label{cor4.5}
The backward shift $B_w$ is weak*-$q$-frequently hypercyclic on
$\ell^\infty$ if $\displaystyle \sum_{n\in\mathbb{N}} \frac {1}
{w_{1}...w_{j+n^q}}$ converges for each $j\in \mathbb{N}$.
\end{corollary}

\begin{proof}
Immediate since $\{e_n:n\geq 1\}$ is a symmetric Schauder basis for
$\ell^1$.
\end{proof}

Thus if $\displaystyle \sum_{n\in\mathbb{N}} \frac {1}
{w_{1}w_{2}...w_{n}}$ converges, $B_w$ is weak*-frequently
hypercyclic on $\ell^\infty$. In fact, the following stronger result
holds.

\begin{proposition}\label{prop4.6}
(1)If $\displaystyle \lim_{n\to \infty}(w_1w_2...w_n)=\infty$, then
the unilateral backward shift $B_w$ is weak*-frequently hypercyclic
on $\ell^\infty$. \\
(2)If $\displaystyle \lim_{n\to \infty}(w_1w_2...w_n)=\infty$ and
$\displaystyle \lim_{n\to \infty}(w_{-1}w_{-2}...w_{-n})=0$, then
the bilateral backward shift $T_w$ is weak*-frequently hypercyclic
on $\ell^\infty(\mathbb{Z})$.
\end{proposition}

\begin{proof}
Since $\displaystyle \lim_{n\to \infty}(w_1w_2...w_n)=\infty$, the
weighted shift $B_w$ is frequently hypercyclic on $c_0$, cf. \cite
{Bayart Grivaux} or \cite {Bayart Matheron}. Also, it is easy to see
that the identity operator from  $c_0$ to $\ell^\infty$ is norm to
weak* continuous and has weak$^*$-dense range. Hence $B_w$ is
weak*-frequently hypercyclic on $\ell^\infty$, by Proposition
\ref{prop2.3}. Similarly, the identity map takes $c_0(\mathbb{Z})$
into $\ell^\infty(\mathbb{Z})$ continuously and densely. Thus $T_w$
is weak*-frequently hypercyclic on $\ell^\infty(\mathbb{Z})$.

\end{proof}

\begin{remark}
We would like to mention here that hypercyclicity on $\ell^\infty$
which is weak$^*$-separable, was studied in \cite {Bes Chan Sanders}
and \cite {Petterson}. It was proved that a backward shift $B_w$ is
weak$^*$-hypercyclic on  $\ell^\infty$ if and only if $\limsup_{n\to
\infty}(w_1w_2...w_n)=\infty$. However, the condition
$\lim_{n\rightarrow \infty} (w_1w_2..w_n)=\infty$ is not necessary
for $B_w$ to be weak$^*$-frequently hypercyclic on $\ell^\infty$.
Indeed, there exists a frequently hypercyclic $B_w$ on $c_0$ (and
thus weak$^*$-frequently hypercyclic on $\ell^\infty$) such that
$w_1w_2..w_n\nrightarrow \infty$. cf. \cite{Bayart Griv}, p. 205.
\end{remark}

\begin{remark}
In view of Theorem \ref{thm3.5}, a weakly $q$-frequently hypercyclic
operator on a separable Banach space, which satisfies the
$q$-frequent hypercyclicity criterion with respect to a weakly dense
set is necessarily norm-$q$-frequently hypercyclic; for the closed
convex sets are the same in the weak and norm topologies .
\end{remark}

Finally, we consider the frequent hypercyclicity of a
non-convolution operator. It is known that any convolution operator
on $H(\mathbb{C})$ (a continuous linear operator that commutes with
all translations) that is not a multiple of the identity operator is
frequently hypercyclic \cite {Bonilla Erdmann1}. The operator
$T_\mu(f)(z)=f'(\mu z)$ on the space $H(\mathbb{C})$ is a
non-convolution for $\mu\not=1$ and was shown to be hypercyclic for
$|\mu|\geq 1$, cf. \cite{Gustavo Hallack} and \cite{Aron Dinesh}. In
fact, this operator is a weighted backward shift with weights
$w_n=n\mu^{n-1}$. So, our result can be derived using the
Proposition \ref{prop4.1}. We rather prove this in the following
way.

\begin{proposition}
Let $H(\mathbb{C})$ be equipped with the compact-open topology. Then
the  operator $T_\mu$ is frequently hypercyclic on $H(\mathbb{C})$
for $|\mu|\geq 1$.
\end{proposition}
\begin{proof}

Let $D$ be the set of all polynomials. Define the map $S_\mu$ by,
\begin{center}
$\displaystyle S_\mu(f)(z)=\mu \int_0^{\frac {1} {\mu} z} f(u)\,du
$.
\end{center}

It is easy to see that $\sum T_\mu ^n(f)$ is absolutely convergent
in $H(\mathbb{C})$ and that $T_\mu S_\mu = I$, the identity on the
set $D$. We fix a $k\geq0$ and consider the function $z^k$. Then
\begin{center}
$\displaystyle S_\mu ^n(z^k)= \frac {k!z^{k+n}} {(k+n)!
\mu^{nk+n(n-1)/2}}.$
\end{center}
Since $|\mu|\geq 1$, the series $\sum S_\mu ^n(f)$ is absolutely
convergent for any polynomial $f$ and we conclude that $T_\mu$ is
frequently hypercyclic by Theorem \ref{thm2.2}.
\end{proof}

We ask the following question. Can one say that the operators
$T_{\mu,b}(f)(z)=f'(\mu z+b)$ on $H(\mathbb{C})$ for $|\mu|\geq 1$
and $b\in \mathbb{C}$, $b\not=0$, are frequently hypercyclic? These
have been shown to be hypercyclic in \cite{Aron Dinesh} and
\cite{Gustavo Hallack}.

\section{Rotation and Powers}

We now remark that rotations and powers of a $q$-frequently
hypercyclic operator on an arbitrary separable topological vector
space remain $q$-frequently hypercyclic for any $q\in \mathbb{N}$.
They also share the same set of $q$-frequently hypercyclic vectors.
These results have been proved for the case $q=1$ in \cite{Bayart
Matheron}. The hypercyclicity of powers and rotations have been
considered by S. I. Ansari \cite{Ansari} and F. Leon-Saavedra and V.
M\"{u}ller \cite{Leon Muller} respectively. For establishing the
following theorem on powers and rotations of $q$-frequently
hypercyclic operators, we need a lemma stated as

\begin{lemma}
Let $A\subset \mathbb{N}$ have positive $q$-lower density and
$\displaystyle \bigcup_{j=1}^{k}I_j=\mathbb{N}$. If $n_1,...,n_k$
are finitely many natural numbers, then
\begin{center}
$\displaystyle \bigcup_{j=1}^{k}(n_j+A\cap I_j)$
\end{center}
has positive $q$-lower density.
\end{lemma}
\begin{proof}
Omitted as it follows on the same lines as given in \cite{Bayart
Matheron}, p. 148.
\end{proof}

Let us denote by $qFHC(T)$, the set of all $q$-frequently
hypercyclic vectors for $T$ and $S^1$, the unit circle in the
complex plane. Then we have

\begin{theorem}
Let $T$ be a $q$-frequently hypercyclic operator on a complex
topological vector space $X$. Then $\lambda T$ and $T^p$ are
$q$-frequently hypercyclic for each $\lambda \in S^1$ and $p\in
\mathbb{N}$. Also $qFHC(T)=qFHC(\lambda T)=qFHC(T^p)$.
\end{theorem}
\begin{proof}
To get this result, we proceed on similar lines as in the case of
frequent hypercyclicity, \cite{Bayart Matheron}, p. 148.

\end{proof}
 Finally, we would like to mention that the notion of $q$-frequent
hypercyclicity is a particular case of $(m_k)$-hypercyclicity
studied in \cite{Bayart Mathero}. Indeed, for a strictly increasing
sequence $(m_k)$ of natural numbers, an element $x\in X$ is called
$(m_k)$-hypercyclic for an operator $T$ on $X$ if for every
non-empty open set $U\subset X$, there exists a strictly increasing
$(n_k)=O(m_k)$ such that $T^{n_k}(x)\in U$ for all $k$. Thus the
case $m_k=k^q$, for all $k$ coincides with the notion of
$q$-frequent hypercyclicity; however, the results in our paper have
no overlap with the results of \cite{Bayart Mathero} except that the
$q$-frequent hypercyclicity $(q\geq 2)$ of the Bergman shift has
been proved in \cite{Bayart Mathero} using the notion of a
hypercyclicity set; see Example 5.3, \cite{Bayart Mathero}.

\textbf{Acknowledgement.}  The authors are thankful to the referee for
his/her careful reading of the paper, and pointing out the
references \cite{Bayart Griv} and \cite{Bayart Mathero}, which
respectively helped them to answer the converse of Proposition
\ref{prop4.6}(1) and mention the $(m_k)$-hypercyclicity. The second
author acknowledges a financial support from the Council of
Scientific and Industrial Research India for carrying out research
at IIT Kanpur.

\bibliographystyle{amsplain}

\end{document}